\documentclass[titlepage,14pt]{extarticle}
\usepackage{amsfonts,amsmath,amssymb,indentfirst}

\usepackage{geometry}
\def\inbr#1{ \{ #1 \} }
\newtheorem{theorem}{Theorem}
\newtheorem{lemma}{Lemma}
\newtheorem{corollary}{Corollary}
\def\proof{{\bf Proof. }}

\allowdisplaybreaks[4]

\begin{document}

\noindent
{\large \bf \scshape On a Certain Integral Over a Triangle} \\
{\bf \scshape Sergey Zlobin}

\vskip 0.5cm

Consider the integral
$$
I_n=\int_{T} \frac{(-\ln xy)^n}{xy},
$$
where T is the triangle with vertices (1,0), (0,1), (1,1).

Its first values are
$$
I_0=\zeta(2), \quad
I_1=2 \zeta(3), \quad
I_2=\frac{9}{2}\zeta(4),
$$
$$
I_3=36 \zeta(5)-12 \zeta(2) \zeta(3), \quad
I_4=\frac{237}{2} \zeta(6) - 48 \zeta(3)^2.
$$
In the general case $I_n$ is equal to a linear combination of
multiple zeta values
$$
\zeta(s_1, s_2, \dots, s_l) = \sum_{n_1 > n_2 > \cdots > n_l \ge
1} \frac{1}{n_1^{s_1} \cdots n_l^{s_l}}
$$
of weight $n+2$ with rational coefficients:

\begin{theorem}
\label{Th1}
For any integer $n \ge 0$ the following identity holds:
$$
I_n=n! \sum_{k=0}^n
\zeta(n-k+2,\inbr{1}_k)
$$
(\; $\inbr{a}_k$ means $\underbrace{a, \dots, a}_k$ \;).
\end{theorem}

\begin{lemma}
\label{lemma1}
Let $k$ and $l$ be integers such that $k \ge 0$, $l \ge 1$. Then the
identity
$$
\int_{0}^1 \frac{(-\ln (1-x))^k}{1-x}   \cdot
(-\ln x)^l dx=k! l! \zeta(l+1,\inbr{1}_k)
$$
holds.
\end{lemma}
\proof Denote the integral in the statement of the lemma by $J$.
Consider the following series expansion:
$$
\frac{(-\ln (1-x))^k}{1-x}=  k!
\sum_{n_1 \ge n_2 > n_3 > \dots > n_{k+1} \ge 1}
\frac{x^{n_1}}{n_2 \cdots n_{k+1}}.
$$
Substituting it into the integral, we get
\begin{align*}
J&=k!
\sum_{n_1 \ge n_2 > n_3 > \dots > n_{k+1} \ge 1}
\frac{1}{n_2 \cdots n_{k+1}}
\int_{0}^1 x^{n_1} (-\ln x)^l dx \\
& = k!
\sum_{n_1 \ge n_2 > n_3 > \dots > n_{k+1} \ge 1}
\frac{1}{n_2 \cdots n_{k+1}} \cdot
\frac{l!}{(n_1+1)^{l+1}}= \\
& = k! l!
\sum_{n_1 > n_2 > n_3 > \dots > n_{k+1} \ge 1}
\frac{1}{n_1^{l+1} n_2 \cdots n_{k+1}} \\
&=k! l! \zeta(l+1,\inbr{1}_k).
\end{align*}

{\bf Remark. }
After obtaining this result, it was found that the integral in the 
statement of the lemma can always be represented
as a polynomial of several variables with rational coefficients in values
of the Riemann zeta function at integers.
The explicit formula is given in \cite{Kolbig1}:
\begin{multline*}
\int_{0}^1 \frac{(-\ln (1-x))^k}{1-x} \cdot
(-\ln x)^l dx \\
=k! l! 
\sum_{p=1}^l \frac{(-1)^{p+1}}{p!}
\sum_{t_i} \frac{H_l(t_1,\dots,t_p)}{t_1 \cdots t_p} 
\zeta(t_1) \cdots \zeta(t_p),
\end{multline*}
where 
$$
H_l(t_1,\dots,t_p)=\sum_{l_i} \binom{t_1}{l_1} \cdots \binom{t_p}{l_p}.
$$
The sum over $t_i$ is to be taken over all sets of integers
$\{ t_i \}$ $(i = 1, \dots ,p)$ which satisfy 
$$
t_i \ge 2, \quad \sum_{i=1}^{p} t_i = k+l+1,
$$
and the sum over $l_i$ over all sets of integers $\{ l_i \}$ $(i = 1, \dots ,p)$
which satisfy 
$$
1 \le l_i \le t_i -1, \quad \sum_{i=1}^p l_i = l.
$$
Our result is shorter and looks more elegant, but the expression involving
only ordinary zeta values is more classical. 
Combining both expressions for the integral we obtain the interesting
\begin{corollary}
If \, $m \ge 2$ and $k \ge 0$, then the value $\zeta(m, \inbr{1}_k )$
can be explicitly represented as a polynomial of several variables 
with rational coefficients in values of the Riemann zeta function at integers.
\end{corollary}
This corollary, including the polynomial formula for
$\zeta(m, \inbr{1}_k )$ (at least implicitly),
has been also proved in \cite{bradley}, using the generating function
$$
\sum_{k, l \ge 0} x^{k+1} y^{l+1} \zeta(l+2, \inbr{1}_k )=
1-\exp \left\{ \sum_{n \ge 2} \frac{x^n+y^n-(x+y)^n}{n} \zeta(n) \right\}.
$$

{\bf Proof} of Theorem \ref{Th1}.
We have
\begin{align*}
I_n&=\sum_{k=0}^n \binom{n}{k} \int_{0}^1 \frac{(-\ln x)^k}{x}
\int_{1-x}^1 \frac{(-\ln y)^{n-k}}{y} dy dx \\
& =\sum_{k=0}^n \binom{n}{k} \int_{0}^1 \frac{(-\ln x)^k}{x} \cdot
\frac{(-\ln (1-x))^{n-k+1}}{n-k+1} dx \\
& =\sum_{k=0}^n \binom{n}{k} \int_{0}^1 \frac{(-\ln
(1-x))^k}{1-x}   \cdot \frac{(-\ln x)^{n-k+1}}{n-k+1} dx.
\end{align*}
Applying Lemma \ref{lemma1}, we get
$$
I_n=\sum_{k=0}^n \binom{n}{k} \cdot \frac{1}{n-k+1} \cdot
k! (n-k+1)! \zeta(n-k+2,\inbr{1}_k).
$$
This is equivalent to the required assertion.

For example, Theorem \ref{Th1} yields
$$
I_2=2(\zeta(4)+\zeta(3,1)+\zeta(2,1,1))=
2\left(\zeta(4)+\frac{\zeta(4)}{4}+\zeta(4)\right)=\frac{9}{2}\zeta(4).
$$
Note that the summands in the sum
$$
\sum_{k=0}^n \zeta(n-k+2,\inbr{1}_k)
$$
have the symmetry
$$
\zeta(n-k+2,\inbr{1}_k)=\zeta(k+2,\inbr{1}_{n-k})
$$
by the duality theorem.

It would be interesting to find a generalization of Theorem
\ref{Th1} with an $m$-dimensional integral for $m \ge 2$.

I thank kindly J. Sondow for raising the problem and several
suggestions.


\newcommand{\namefont}{\scshape}
\newcommand{\titlefont}{\itshape}
\def\nomer{No.}

\end {document}